\magnification=1200
\def\defeq{\buildrel \rm def \over =}
\def\Scr{\cal}
\def\tpof{$({\bf 3+1})$-free}

\catcode`\@=11
\font\teneuf=eufm10
\font\seveneuf=eufm7
\font\fiveeuf=eufm5
\newfam\euffam
\textfont\euffam=\teneuf \scriptfont\euffam=\seveneuf
  \scriptscriptfont\euffam=\fiveeuf
\def\euf@{\hexnumber@\euffam}
\def\Frak{\relax\ifmmode\let\next\Frak@\else
 \def\next{\errmessage{Use \string\Frak\space only in math mode}}\fi\next}
\def\Frak@#1{{\Frak@@{#1}}}
\def\Frak@@#1{\fam\euffam#1}
\catcode`\@=12

\newcount\thmno \thmno=1
\long\def\theorem#1#2{\edef#1{Theorem~\the\thmno}\noindent
    {\sc Theorem \the\thmno.\enspace}{\it #2}\endgraf
    \penalty55\global\advance\thmno by 1}
\newcount\lemmano \lemmano=1
\long\def\lemma#1#2{\edef#1{Lemma~\the\lemmano}\noindent
    {\sc Lemma \the\lemmano.\enspace}{\it #2}\endgraf
    \penalty55\global\advance\lemmano by 1}
\newcount\propno \propno=1
\long\def\proposition#1#2{\edef#1{Proposition~\the\propno}\noindent
    {\sc Proposition \the\propno.\enspace}{\it #2}\endgraf
    \penalty55\global\advance\propno by 1}
\newcount\corolno \corolno=1
\long\def\corollary#1#2{\edef#1{Corollary~\the\corolno}\noindent
    {\sc Corollary \the\corolno.\enspace}{\it #2}\endgraf
    \penalty55\global\advance\corolno by 1}
\newcount\conjno \conjno=1
\long\def\conjecture#1#2{\edef#1{Conjecture~\the\conjno}\noindent
    {\sc Conjecture \the\conjno.\enspace}{\it #2}\endgraf
    \global\advance\conjno by 1}

\newcount\eqnno \eqnno=0
\def\eqn#1{\global\edef#1{(\the\secno.\the\eqnno)}#1\global
   \advance\eqnno by 1}

\newcount\secno \secno=0
\def\sec#1\par{\eqnno=1\global\advance\secno by 1\bigskip
    \noindent\centerline{\sc\the\secno. #1}\par\nobreak\noindent}

\def\definition{\noindent{\sc Definition.\enspace}}
\def\proof{\noindent{\it Proof.\enspace}}
\def\qed{~\hbox{\quad\vbox{\hrule height.4pt
    \hbox{\vrule width.4pt height7pt \kern7pt \vrule width.4pt}
    \hrule height.4pt}}}

\font\noterm=cmr9

\catcode`\@=11
\def\footnote#1{\edef\@sf{\spacefactor\the\spacefactor}#1\@sf
      \insert\footins\bgroup\noterm
      \interlinepenalty100 \let\par=\endgraf
        \leftskip=\z@skip \rightskip=\z@skip
        \splittopskip=10pt plus 1pt minus 1pt \floatingpenalty=20000
        \medskip\item{#1}\bgroup\strut\aftergroup\@foot\let\next}
\catcode`@=12
\dimen\footins=30pc 
\newdimen\notespc \notespc=11pt
\newcount\notenum \notenum=0
\newdimen\currentspc
\def\note#1{\currentspc=\the\baselineskip
      \parskip=0pt\baselineskip=\notespc
      \global\advance\notenum by 1\footnote{$^{\the\notenum}$}\bgroup
      #1\egroup { }\baselineskip=\currentspc \parskip=12pt plus4pt minus4pt}

\parskip=12pt plus4pt minus4pt

\font\sc=cmcsc10
\def\itc#1{{\it #1\/}}
\def\X{{\rm X}}


\vbox to \vsize
{
\vfill
\centerline{\bf DESCENTS, QUASI-SYMMETRIC FUNCTIONS,}
\centerline{\bf AND THE CHROMATIC SYMMETRIC FUNCTION}
\bigskip\bigskip
\centerline{Timothy Y. Chow}
\centerline{Department of Mathematics}
\centerline{University of Michigan at Ann Arbor}
\vfill\vfill
}
\eject

\sec Introduction

The theory of $P$-partitions continues to spawn new ideas more than
twenty years after its birth.  Our main object of interest here is
one such outgrowth, namely the expansion of Stanley's chromatic
symmetric function in terms of Gessel's fundamental quasi-symmetric
functions~$Q_{S,d}$ (reproduced as Theorem~1 below).
Although innocent-looking,
this expansion has numerous ramifications,
some of them surprising.
The purpose of this paper is to explore some of these offshoots.

In section~3, we recall the result,
stating it in a way that differs slightly from the usual formulation;
the justification for this modification of standard terminology is that
it shows more clearly the relationship with two other closely related
results in the literature:
Chung and Graham's $G$-descent expansion of the chromatic polynomial
[3, Theorem~2]
and the expansion of the path-cycle symmetric function
in terms of the~$Q_{S,d}$ [2, Proposition~7].
The original proofs of these latter two results did not appeal directly
to Stanley's expansion; here we show that the $G$-descent result
and an important special case of the path-cycle symmetric function result
are essentially special cases of Stanley's result.
In section~4, we investigate the implications of Theorem~1 for
Robinson-Schensted algorithms for \tpof\ posets,
a topic that has attracted some recent attention
([9, section~3.7] and~[14]).
Finally, in section~5, we investigate the connection with the
new symmetric function basis that was introduced in~[2].

\sec Preliminaries

We shall assume that reader is familiar with the basic facts about
set partitions, posets, permutations, and so on; a good reference is~[12].
Our notation for symmetric functions and partitions for the most
part follows that of Macdonald~[8].
If $\lambda$ is an integer partition, we write $r_\lambda!$
for $r_1!r_2!\cdots$, where $r_i$ is the number of parts of~$\lambda$
of size~$i$.
We will always take our symmetric functions in countably many variables.
In addition to the usual symmetric function bases,
we shall need the
\itc{augmented monomial symmetric functions~$\tilde m_\lambda$}~[4],
which are defined by
$$\tilde m_\lambda \defeq r_\lambda!\, m_\lambda,$$
where $m_\lambda$ of course denotes the usual monomial symmetric function.
We will sometimes use
\itc{set} partitions instead of \itc{integer} partitions in subscripts;
for example, if $\pi$ is a set partition then the expression $p_\pi$ is
to be understood as an abbreviation for~$p_{\rm type(\pi)}$.
We will use $\omega$ to denote the involution that sends $s_\lambda$
to~$s_{\lambda'}$.

Throughout, the unadorned term \itc{graph} will mean a finite simple
labelled undirected graph.  If $G$ is a graph we let $V(G)$ and~$E(G)$
denote its vertex set and edge set respectively.
A \itc{stable partition} of~$G$ is a partition of~$V(G)$ such that
every block is a stable set, i.e.,
no two vertices in the same block are connected by an edge.
Stanley's \itc{chromatic symmetric function~$\X_G$} is defined by
$$\X_G \defeq \sum_\pi \tilde m_\pi,$$
where the sum is over all stable partitions~$\pi$ of~$G$.
For motivation for the definition of~$\X_G$, see~[11].

If $d$ is a positive integer, we use the notation
$[d]$ for the set $\{1,2,\ldots,d\}$.

Following Gessel~[7] and Stanley~[11], we define a power series
in the countably many variables $x=\{x_1,x_2,\ldots\}$
to be \itc{quasi-symmetric} if the coefficients of
$$x_{i_1}^{r_1}x_{i_2}^{r_2}\cdots x_{i_k}^{r_k}
  \qquad{\rm and}\qquad
  x_{j_1}^{r_1}x_{j_2}^{r_2}\cdots x_{j_k}^{r_k}$$
are equal whenever $i_1<i_2<\cdots<i_k$ and $j_1<j_2<\cdots<j_k$.
For any subset~$S$ of~$[d-1]$ define the \itc{fundamental}
quasi-symmetric function $Q_{S,d}(x)$ by
$$Q_{S,d}(x) = \sum_{\scriptstyle i_1\le\cdots\le i_d \atop
  \scriptstyle i_j < i_{j+1}\;{\rm if}\; j\in S}
  x_{i_1} x_{i_2} \cdots x_{i_d}.$$
Sometimes we will write $Q_{S,d}$ for $Q_{S,d}(x)$
if there is no danger of confusion.

If $g$ is a symmetric or quasi-symmetric function
in countably many variables and of bounded degree,
then we shall write $g(1^n)$ for the polynomial in the variable~$n$
obtained by setting $n$ of the variables equal to one
and the rest equal to zero.
An important example of this procedure is given in the following
proposition, whose (easy) proof we leave as an exercise.

\proposition\specialQ{For any $S\subseteq [d-1]$,
$$Q_{S,d}(1^n) = {n + d - |S| - 1 \choose d}.$$}

\def\fo{{\Frak o}}
\sec The Fundamental Theorem

The fundamental result in this subject is Stanley's expansion of
the chromatic symmetric function in terms of Gessel's fundamental
quasi-symmetric functions.
We shall now present this result;
more precisely, as mentioned in the introduction,
we shall present a reformulation of the result,
and then we will go on to show how this reformulation subsumes
Chung and Graham's $G$-descent expansion of the chromatic polynomial
and a special case of the expansion of the path-cycle symmetric function
in terms of the~$Q_{S,d}$.

A number of details will be omitted from the proofs in this section
because the arguments consist mostly of definition-chasing.

We need some definitions.
The first of these looks trivial but is actually one of the most important.

\definition
A \itc{sequencing} of a graph or a poset with a vertex set~$V$ that
has cardinality~$d$ is a bijection $s:[d]\to V$.

It is helpful to think of a
sequencing as the sequence $s(1), s(2), \ldots, s(d)$ of vertices.
The reason we claim that this definition is important is that the
usual approach to this subject regards a permutation of some kind
(either of~$[d]$ or of~$V$) as the fundamental object of interest,
but as we shall see below, it is often sequencings that are most
natural to consider.  Even the standard approach often finds it
necessary to resort to inverse maps at certain points to convert
permutations to sequencings; by focusing on sequencings directly
we obviate this.

``Dual'' to the notion of a sequencing is a
\itc{labelling,} which is a bijection $\alpha:V\to [d]$.
A labelling~$\alpha$ of a poset is \itc{order-reversing} if
$\alpha(x)>\alpha(y)$ whenever $x<y$.

A sequencing~$s$ of a graph~$G$ induces an acyclic orientation~$\fo(s)$
of~$G$: if $i<j$ and $s(i)$ is adjacent to~$s(j)$, then direct the
edge from $s(j)$ to~$s(i)$.  The acyclic orientation in turn induces
a poset structure~$\overline{\fo}(s)$ on the vertex set of~$G$:
make $s(i)$ less than $s(j)$ whenever $s(j)$ points to~$s(i)$
and then take the transitive closure of this relation.

Let $G$ be a graph with $d$ vertices.
If $\alpha$ is a labelling of~$G$ and $s$ is a sequencing of~$G$,
then we say that
\itc{$s$ has an $\alpha$-descent at~$i$} (for $i\in [d-1]$)
if the permutation $\alpha \circ s$ has a descent at~$i$.
The \itc{$\alpha$-descent set $D(\alpha,s)$ of~$s$} is the set
$$\{i\in [d-1] \mid \hbox{$s$ has an $\alpha$-descent at~$i$}\}.$$
(It is helpful to visualize this
by visualizing a numerical label
on each element of the sequence $s(1), s(2), \ldots, s(d)$;
the sequence of labels is
the one-line representation of the permutation $\alpha\circ s$
and the descents occur at the descents of this permutation.)

We can now state Stanley's theorem.

\theorem\fund{Let $G$ be a graph with $d$ vertices.
Suppose that to each sequencing~$s$ of~$G$ there is associated an
order-reversing labelling~$\alpha_s$ of~$\overline{\fo}(s)$.
Suppose further that $\alpha_s = \alpha_{s'}$
whenever $s$ and~$s'$ are two sequencings of~$G$
that induce the same acyclic orientation of~$G$.  Then
$$\X_G = \sum_{{\rm all\ sequencings}\ s} Q_{D(\alpha_s,s),d}.$$}

\noindent{\it Sketch of proof.\enspace}
For each acyclic orientation~$\fo$ of~$G$,
let $s$ be some sequencing that induces~$\fo$ and define
$\omega_\fo : \overline{\fo} \to [d]$
to be the order-reversing bijection~$\alpha_s$.
Let ${\Scr L}(\overline{\fo},\omega_\fo)$ be
the set of all linear extensions of~$\overline{\fo}$,
regarded as permutations of~$[d]$ via~$\omega_\fo$,
and if $e$ is a permutation let $D(e)$ denote the descent set of~$e$.
Then [11, Theorem~3.1], combined with [11, equation~(8)], states that
$$\eqalignno{\X_G &= \sum_{\fo}
   \sum_{e \in {\Scr L}(\overline{\fo},\omega_\fo)}
    Q_{D(e),d},&\eqn\fundeqn\cr}$$
where the first sum is over all acyclic orientations~$\fo$ of~$G$.

Now there is a bijection between the set of all sequencings
of~$G$ and the set of ordered pairs
$\bigl\{(\fo, e) \mid e \in{\Scr L}(\overline{\fo})\bigr\}$---given
a sequencing~$s$,
let $\fo$ be the acyclic orientation induced by~$s$ and let
$e = \omega_\fo \circ s = \alpha_s \circ s$.
\fund\ then follows from~\fundeqn once we verify that
$D(e)$ corresponds to $D(\alpha_s,s)$ under this bijection.\qed

Chung and Graham~[3, Theorem~2]
have shown that when the chromatic polynomial of a graph is expanded
in terms of the polynomial basis
$${x+k\choose d}_{k=0,\ldots,d},$$
then the coefficients can be interpreted in terms of what they call
\itc{$G$-descents.}
In their paper, Chung and Graham give a sketch of
a somewhat complicated proof of this result,
and remark that while in principle it follows from
Brenti's expansion~[1, Theorem~4.4]
(which in turn is essentially what one obtains by specializing \fund\ 
via the map $g\mapsto g(1^n)$),
the implication is not particularly direct.
However, Chung and Graham's result follows directly from \fund\
by choosing the $\alpha_s$ appropriately and then
specializing from symmetric functions to one-variable polynomials,
as we shall now see.

Again, we need some definitions.
To \itc{peel} a poset~$P$ is to remove its minimal elements,
then to remove the minimal elements of what is left, and so on.
The \itc{rank~$\rho(x)$} of an element $x\in P$ is the stage
at which it is removed in the peeling process.

Next we give the definition of Chung and Graham's concept of a $G$-descent,
translated into our terminology.
Let $G$ be a graph with $d$ vertices.
Let $\beta$ be a labelling of~$G$ and
let $s$ be a sequencing of~$G$.
If $v$ is a vertex of~$G$ then we define $\rho(v)$
by using the poset structure~$\overline{\fo}(s)$.
We then say that \itc{$s$ has a
CG~$\beta$-ascent at~$i$} (for $i\in [d-1]$) if either

\item{1.} $\rho(s(i)) < \rho(s(i+1))$ or

\item{2.} $\rho(s(i)) = \rho(s(i+1))$ and $\beta(s(i)) < \beta(s(i+1))$.

\noindent
The \itc{CG~$\beta$-ascent set of~$s$} is defined in the obvious way.
We then have the following result.

\corollary\CG{If $G$ is a graph with $d$ vertices and a labelling~$\beta$,
then
$$\eqalignno{\X_G &= \sum_S N_S\, Q_{S,d},&\eqn\CGeqn\cr}$$
where the sum is over all subsets $S\subseteq [d-1]$ and $N_S$ is the number
of sequencings of~$G$ with CG~$\beta$-ascent set~$S$.}

\noindent{\it Sketch of proof.\enspace}
The appropriate choices of~$\alpha_s$ in \fund\ are as follows.
Given a sequencing~$s$, arrange the vertices of~$G$ in the
following ``peeling order'':
first take the elements of highest rank in~$\overline{\fo}(s)$,
then the elements of next highest rank, and so on;
arrange elements with the same rank
in decreasing order of their $\beta$-labels.
Now define the labelling~$\alpha_s$ by setting $\alpha_s(v) = j$ 
where $j$ is the position of~$v$ in the peeling order.
It is now straightforward to check that the CG $\beta$-ascent set of~$s$
coincides with the $\alpha_s$-descent set of~$s$.\qed

Chung and Graham's result [3, Theorem~2]
now follows as a special case of \CG.
For if we apply the map $g\mapsto g(1^n)$ to \CGeqn,
the left-hand side specializes to the chromatic polynomial of~$G$
([11, Proposition~2.2]) and by \specialQ\ the right-hand side specializes
to the binomial coefficient sum
$$\sum_S N_S {n+d-|S|-1 \choose d} = \sum_k N_k {n+k \choose d},$$
where $N_k$ is the number of sequencings with $d-1-k$ CG~$\beta$-ascents,
i.e., with $k$ CG~$\beta$-descents (where CG~$\beta$-descents are 
defined in the natural way).  This is exactly [3, Theorem~2].

Our second corollary involves the expansion of the
path-cycle symmetric function~$\Xi_D$
in terms of the~$Q_{S,d}$ [2, Proposition~7].
(We shall not give the formal definition of the path-cycle symmetric
function here because we will not need it;
suffice it to say that it is a certain symmetric function
invariant~$\Xi_D$ that can be associated to any digraph~$D$.)
For certain digraphs~$D$, $\Xi_D$ coincides with
the chromatic symmetric function~$\X_G$ of some graph~$G$,
and therefore, in these cases,
[2, Proposition~7] gives an interpretation of the coefficients
of the $Q_{S,d}$-expansion of~$\X_G$.
This interpretation is ostensibly different from the one given by \fund,
but as we shall show presently, it again follows directly from \fund\
via suitable choices of~$\alpha_s$.

We shall now make these somewhat vague remarks precise.
Let $P$ be a poset with $d$ vertices.
If $s$ is a sequencing of~$P$,
we say that $s$ has a \itc{descent at~$i$} (for $i\in [d-1]$)
if $s(i)\not< s(i+1)$.  The \itc{descent set~$D(s)$ of~$s$}
is again defined in the obvious way.
The \itc{incomparability graph~${\rm inc}(P)$} of~$P$ is the
graph with the same vertex set as~$P$ and in which two vertices
are adjacent if and only if they are incomparable elements of~$P$.

An acyclic, transitively closed digraph is equivalent to a poset.
According to [2, Proposition~2], the path-cycle symmetric function
of such a digraph coincides with the chromatic symmetric function
of the incomparability graph of the equivalent poset.
Therefore, what [2, Proposition~7] says in this case is the following.

\corollary\pathcyc{Let $P$ be a poset with $d$ vertices.  Then
$$\X_{{\rm inc}(P)} =
    \sum_{{\rm all\ sequencings}\ s} Q_{D(s),d}.$$}

We now claim that this result can also be derived from \fund.

\noindent{\it Sketch of proof.\enspace}
We define the $\alpha_s$ as follows.  Let $s$ be any sequencing
of~${\rm inc}(P)$.
The maximal elements of~$\overline{\fo}(s)$
form a stable set in~${\rm inc}(P)$ and therefore a chain in~$P$;
call the minimal (with respect to the ordering of~$P$,
not of~$\overline{\fo}(s)$) element of this chain~$v_1$, and
set $\alpha_s(v_1)=1$.  Now delete~$v_1$ and repeat the procedure,
i.e., let $v_2$ be the $P$-minimal element
among the $\overline{\fo}(s)$-maximal elements of the deleted graph,
and set $\alpha_s(v_2)=2$.
Continue in this way until $\alpha_s(v)$ is defined for all~$v$.
We leave to the reader the
(straightforward although not entirely trivial)
task of verifying that the $\alpha_s$-descents of~$s$
(considered as a sequencing of~${\rm inc}(P)$)
coincide with the descents of~$s$
(considered as a sequencing of~$P$).\qed

\sec Robinson-Schensted and \tpof\ Posets

A poset is said to be \itc{\tpof} if it contains no induced
subposet isomorphic to the disjoint union of a three-element chain with
a singleton.  Unless otherwise noted, all posets in this section are
assumed to be \tpof.

Gasharov~[6] has proved a remarkable result about the
expansion of $\X_{{\rm inc}(P)}$ in terms of Schur functions.
To state it we must first recall the notion of a
$P$-tableau.  If $P$ is any poset, a \itc{(standard) $P$-tableau}
is an arrangement of the elements of~$P$ into a Ferrers shape
such that the rows are strictly increasing (i.e., each row is a chain)
and the columns are weakly increasing
(by which we mean that if $u$ appears \itc{immediately above~$v$}
[in English notation] then $u\not>v$).
Each element of~$P$ appears exactly once in the tableau.
Then Gasharov's result is the following.

\theorem\gasharov{Let $P$ be a \tpof\ poset.  Then
$$\X_{{\rm inc}(P)} = \sum_\lambda f^\lambda_P s_\lambda,$$
where $f^\lambda_P$ is the number of $P$-tableaux of shape~$\lambda$.}

It would be nice to have a direct bijective proof of \gasharov\
(Gasharov's proof is not).
In~[13] Stanley remarks that
when $P$ is a chain, $f^\lambda_P$ is just the number of
standard Young tableaux of shape~$\lambda$, so a bijective proof of
\gasharov\ is provided by the Robinson-Schensted correspondence.
(For background on Robinson-Schensted and tableaux, see~[10].)
Stanley further remarks that Magid~[9, Section~3.7]
has produced a generalization
of the Robinson-Schensted correspondence that provides the desired
bijective proof of \gasharov.
However, the exposition in~[9, Section~3.7] is difficult to follow,
and to the best of my understanding there is an error in the construction.
Let $P$ be the four-element poset whose Hasse diagram looks
like an uppercase~``N.''  (We shall refer to this poset as \itc{Poset~N.})
Label the vertices $a$, $b$, $c$, and~$d$ from left to right, top to bottom
(like reading English).
Then the two sequences $dacb$ and $dbca$ appear to generate the same
pair of tableaux under Magid's insertion algorithm,
which should not happen since the insertion algorithm is supposed
to give a bijection between sequencings of the poset and pairs of tableaux.
It is possible that I am misinterpreting Magid's algorithm,
but if I am correct then the problem of finding a bijective proof
of \gasharov\ is still open.
The best partial result is due to
Sundquist, Wagner, and West~[14], who provide an algorithm that
gives the desired bijection for a certain
proper subclass of \tpof\ posets.

We shall say more about the algorithm in~[14] in a moment,
but our main purpose here is
to observe that combining \pathcyc\ with \gasharov\
gives us some insight into the kind of Robinson-Schensted algorithm
we want.  If $\lambda\vdash d$ then from  [11, equation~(15)] we have
$$s_\lambda = \sum_S f^\lambda_S Q_{S,d},$$
where the sum is over all $S\subseteq [d-1]$ and $f^\lambda_S$ is
the number of standard Young tableaux with shape~$\lambda$
and descent set~$S$.  Combining this with \gasharov\ yields
$$\X_{{\rm inc}(P)} = \sum_S \sum_\lambda f^\lambda_P f^\lambda_S Q_{S,d}.$$
In other words,
the coefficient of~$Q_{S,d}$ in $\X_{{\rm inc}(P)}$
is the number of ordered pairs $(T,T')$ where $T$
is a $P$-tableau and $T'$ is a standard Young tableau with the same shape
and with descent set~$S$.

Comparing this with \pathcyc, we see that not only does there exist a
bijection between sequencings of~$P$ and ordered pairs $(T,T')$
with $T$ a $P$-tableau and $T'$ a standard Young tableau,
but there exists such a bijection with the further property
that it respects descents.
(It is well known that this is true
in the case of the usual Robinson-Schensted algorithm.)
It is therefore natural to hope for an algorithm
that also respects descents.  For one thing, this
would provide an alternative proof of Gasharov's theorem.

We might ask if the Sundquist-Wagner-West algorithm respects descents,
at least for the class of \tpof\ posets to which it is applicable.
The answer is no, and again the abovementioned sequencings of Poset~N
furnish counterexamples.  However, we do have the following result.

\theorem\sww{The Sundquist-Wagner-West algorithm respects descents
when restricted to the class of \tpof\ posets that do not contain
Poset~N as an induced subposet.}

\proof
The Sundquist-Wagner-West algorithm applies to a more general class
of objects than we have been discussing here, but in our present
context, it reduces to the following.  Let $P$ be a \tpof\ poset.
Given a sequencing~$s$ of~$P$, construct an ordered pair $(T,T')$
where $T$ is a $P$-tableau and $T'$ is a standard Young tableau
by \itc{inserting} $s(1)$, $s(2)$, and so on in turn.
The $P$-tableau~$T$ will be the insertion tableau, and~$T'$
will be the recording tableau.
The recording is done in the normal way and requires no comment.
To insert an element~$s(i)$ into~$T$, observe that each row~$R$ of~$T$
is a chain.  (This property is trivial to begin with and it will
be easy to see that it is preserved at each stage of the insertion process.)
Therefore, since $P$ is \tpof, $s(i)$ is incomparable to at most two
elements of~$R$.  If $s(i)$ is incomparable to zero elements of~$R$,
then the situation is indistinguishable from standard Robinson-Schensted,
so proceed in the expected way: append $s(i)$ to the end of~$R$ if
$s(i)$ is greater than every element of~$R$; otherwise, let $s(i)$ bump
the smallest element of~$R$ greater than~$s(i)$ and proceed inductively
by inserting the bumped element into the next row of the insertion tableau.
If $s(i)$ is incomparable to exactly one element of~$R$, make $s(i)$ bump
that one element.  Finally, if $s(i)$ is incomparable to two elements of~$R$,
make $s(i)$ ``skip over''~$R$ and inductively insert $s(i)$ into the next row.
We remark that it is easy to show that if there are any elements in~$R$
incomparable to~$s(i)$, then these elements must be
in a single consecutive block and
that $s(i)$ must be greater than everything to the left of this block
and less than everything to the right of this block.
Keeping this fact in mind will make it easier to follow the arguments
below.

Sundquist, Wagner and West prove that
the above algorithm produces a bijection
if $P$ is what they call ``beast-free'' in addition to being \tpof.
Since the beast contains Poset~N as an induced subposet,
the bijection is valid for the posets that we are concerned with here.

To show that descents are respected in this algorithm,
we proceed by a straightforward case-by-case analysis.
Suppose first that $s(i)\not<s(i+1)$.
We wish to show that $i+1$ appears in a lower row than~$i$ in the recording
tableau.
We claim first that when $s(i+1)$ is inserted,
it cannot be appended at the end of row~1.
To see this, back up and think about
what could have happened when $s(i)$ was inserted.
If $s(i)$ did not skip over row~1, then $s(i+1)$ could not then be
appended to row~1 since $s(i)\not<s(i+1)$.
If on the other hand $s(i)$ did skip over row~1,
then $s(i)$ must be incomparable to two elements in row~1,
and because $P$ is \tpof, we must have $s(i)>s(i+1)$,
and $s(i+1)$ cannot be appended to row~1 because this would force
$s(i)$ to be greater than everything in row~1, contradiction.

Now if $s(i)$ is appended to the end of row~1 then we are done.
Otherwise, each of $s(i)$ and $s(i+1)$ gives rise to an element
to be inserted into row~2;
call these two elements $u$ and~$v$ respectively.
(They need not be distinct from $s(i)$ and~$s(i+1)$ but they must be
distinct from each other.)
By induction it suffices to show that $u\not<v$.  We have
several cases.

\item{1.} Suppose $u=s(i)$, i.e., suppose $s(i)$ skips over row~1.
If $v=s(i+1)$ then we are done.  Otherwise, suppose towards a contradiction
that $s(i)<v$.  Consider the situation before the insertion of~$s(i)$.
Since $v\ne s(i+1)$, $v$ must be in row~1, and since $s(i)<v$, $s(i)$
is less than everything to the right of~$v$.  But $s(i)$ is incomparable
to two elements in row~1, so there must exist at least two elements in
row~1 to the left of~$v$.  Let $q$ and~$r$ be the two elements in row~1
immediately preceding~$v$.  Now $s(i+1)$ bumped~$v$ so $q<r<s(i+1)$.
Since $s(i)\not<s(i+1)$, we have $s(i)\not<q$ and $s(i)\not<r$.
But since $s(i)$ is less than $v$ and everything to the right of~$v$,
we must have $s(i)\not>q$ and $s(i)\not>r$ for otherwise there could
not be two elements in row~1 incomparable to~$s(i)$.  Therefore
$s(i)\not>s(i+1)$ and $s(i)$ together with $q<r<s(i+1)$ is a
($\bf 3+1$), contradiction.

\item{2.} Suppose $u\ne s(i)$ and $s(i)<u$.  If $s(i+1)=v$ then
since $s(i)\not<s(i+1)=v$ and $s(i)<u$ we must have $u\not<v$
and we are done.  So we may assume that $s(i+1)\ne v$.  Suppose
towards a contradiction that $u<v$.  Then when $s(i+1)$ is inserted
into row~1 it bumps something (namely~$v$) that is greater than~$u$
and thus greater than~$s(i)$.
Since $s(i)$ is sitting in row~1 when $s(i+1)$ is inserted,
this forces $s(i+1)>s(i)$, contradiction.

\item{3.} Suppose that $u\ne s(i)$ and that $s(i)$ and~$u$ are incomparable.
We have two subcases: either $s(i+1) \ne v$ or $s(i+1)=v$.
In the former case, suppose towards a contradiction that $u<v$.
Since $v>u$, $v$ must be sitting in row~1 to the right of~$s(i)$
just before $s(i+1)$ bumps it.  Therefore $v>s(i)$ and hence $s(i+1)>s(i)$
(since $s(i+1)$ bumps $v$ and not~$s(i)$), contradiction.
In the latter case, again suppose towards a contradiction
that $u<v$.  Since $s(i)$ and~$u$ are incomparable, this implies that
$s(i+1)=v\not<s(i)$, i.e., that $s(i)$ and $s(i+1)$ are incomparable.
Then $u<s(i+1)$ together with~$s(i)$ form a ($\bf 2+1$), so that $s(i)$
is one of the two elements in row~1 incomparable to~$s(i+1)$ that
cause $s(i+1)$ to skip over row~1.  Let $w$ be the other element in row~1
incomparable to~$s(i+1)$; then $w$ is either the immediate successor
or the immediate predecessor of~$s(i)$---and therefore of~$u$ before
$u$ was bumped by~$s(i)$.  Actually, though, $w$ cannot be a
predecessor of~$u$ since this would make $s(i+1)>w$.
Combining all this information, we see that
$s(i) < w > u < s(i+1)$
together form an induced subposet isomorphic to Poset~N, contradiction.

\noindent
To complete the proof of the theorem we just need to show that if
$s(i)<s(i+1)$ then we do not obtain a descent in the recording tableau.
If $s(i+1)$ is appended to the end of row~1 then we are done.
If $s(i)$ is appended to the end of row~1 then so is~$s(i+1)$
and again we are done.  Therefore, as before, it is enough by
induction to show that ``$u<v$.''

Suppose first that $s(i)$ bumps some element~$u$ from row~1.
If $s(i+1)$ also bumps some element~$v$ from row~1 then
since $s(i)<s(i+1)$ we must have $u<v$, so by induction we are done;
therefore we may assume that $s(i+1)$ skips over row~1.
Suppose towards a contradiction that $u\not<s(i+1)$.
Then $s(i+1)$ is not greater than the element in row~1 immediately
to the right of~$u$, but $s(i+1)$ \itc{is} greater than~$s(i)$,
which displaces~$u$.  Therefore, after the insertion of~$s(i)$,
the two elements in row~1 incomparable to $s(i+1)$ must be the
two elements $q$ and~$r$ in row~1 immediately to the right of~$s(i)$.
Hence $y\not<u$, but then $u<q<r$ and~$y$ form a $({\bf 3+1})$,
contradiction.

It remains to consider the case when $s(i)$ skips over row~1.
If $s(i+1)$ also skips over then we are done.
We have two remaining subcases: either the element~$v$ that $s(i+1)$
bumps is larger than~$s(i+1)$ or else $v$ and $s(i+1)$ are incomparable.

In the former case, let $w_1$ and~$w_2$ be the two elements in row~1
incomparable to~$s(i)$ (just prior to the insertion of~$s(i)$).
Since $s(i)<s(i+1)$ we must have $s(i+1)\not<w_1$ and $s(i+1)\not<w_2$.
Since by assumption $s(i+1)$ bumps something larger than itself,
we must have $y>w_1$ and $y>w_2$.  Therefore $v$ must lie to the
right of $w_1$ and~$w_2$, so $v>s(i)$, which is want we want to show.

In the latter case, suppose towards a contradiction that $s(i)\not<v$.
We cannot have $v<s(i)$ because then $v<s(i)<s(i+1)$, contradicting
the incomparability of~$v$ and $s(i+1)$.  So $v$ is incomparable to~$s(i)$.
Consider row~1 just before the insertion of $s(i+1)$; $s(i)$ is
incomparable to two elements in row~1, and one of these is~$v$.
The other one, which we may call~$w$, must be either the immediate
predecessor or the immediate successor of~$v$.
If $w$ is the immediate successor of~$v$ then this forces $s(i+1)<w$
and since $s(i)<s(i+1)$ this implies $s(i)<w$, contradiction.
Therefore $w$ is the immediate predecessor of~$v$.
Combining this information we see that $s(i)<s(i+1)>w<v$
is an induced subposet isomorphic to Poset~N, contradiction.\qed

Possibly, then, the Sundquist-Wagner-West algorithm needs to be modified
not only in the case of posets containing the ``beast'' but also
beast-free posets that contain Poset~N.  However, so far we have
not been able to find a modification of the Sundquist-Wagner-West
algorithm with all the properties we would like it to have.

\sec The Symmetric Function Basis $\{\xi_\lambda\}$

In~[2] a new symmetric function basis, which we shall denote
by~$\{\xi_\lambda\}$ (in place of the original but more cumbersome
notation~$\{\tilde\Xi_\lambda\}$), is introduced.
For completeness we repeat the definition here.
For each integer partition~$\lambda$,
let $D_\lambda$ denote the digraph
consisting of a disjoint union of directed paths
such that the $i$th directed path has $\lambda_i$ vertices.
If $F$ is a subset of the set~$E(D_\lambda)$ of edges of~$D_\lambda$,
then the spanning subgraph of~$D_\lambda$ with edge set~$F$ is a
disjoint union of directed paths.
The multiset of sizes of these directed paths forms
an integer partition which we denote by $\pi(F)$.
The number of parts of~$\pi(F)$ is denoted by $\ell\bigl(\pi(F)\bigr)$.
Then the symmetric function~$\xi_\lambda$ is defined by
$$\xi_\lambda = \sum_{F\subseteq E(D_\lambda)}
   {\tilde m_{\pi(F)} \over \ell\bigl(\pi(F)\bigr)!},$$
where the sum is over all subsets $F$ of $E(D_\lambda)$.

In [2, Theorem~3] it is stated that $\X_G$ is $\xi$-positive
(i.e., that its expansion
in terms of the $\xi_\lambda$ has nonnegative coefficients).
The proof, however, is not given there.
My original proof of this claim was a direct argument
giving a combinatorial interpretation
of the coefficients in this expansion
in terms of Chung and Graham's $G$-descents.
However, a different proof will be presented here that is perhaps
more illuminating, since it shows how the result follows from \fund.

We need a technical lemma.
If $\pi$ and~$\sigma$ are set partitions,
write $\pi \le \sigma$ for ``$\pi$ refines~$\sigma$.''
If $\pi\le\sigma$, let $k_i$ denote the number of blocks of~$\sigma$
the are composed of $i$~blocks of~$\pi$,
and following Doubilet~[5] define
$$\lambda(\pi,\sigma)! \defeq \prod_i i!^{k_i}.$$
Also, given any integer partitions $\mu$ and~$\nu$,
let $\pi$ be any set partition of type~$\mu$ and define
$$c_{\mu,\nu} \defeq
  \sum_{\{\sigma\ge\pi \mid {\rm type}(\sigma)=\nu\}}
  \lambda(\pi,\sigma)!.$$

\lemma\technical{The number of subsets~$F$ of $E(D_\lambda)$
such that $\pi(F)=\nu$ equals $c_{\nu,\lambda} r_\lambda!/r_\nu!$.}

\proof
See the proof of [2, Proposition~13].\qed

If $S$ is a subset of $[d-1]$ then we define the \itc{type} of~$S$
to be the integer partition whose parts are the lengths of the subwords
obtained by breaking the word $123\ldots d$ after each element of~$S$.

\theorem\xiquasi{Let $g$ be any symmetric function.  If $a_\lambda$
and $b_{S,d}$ are constants such that
$$g = \sum_\lambda a_\lambda \xi_\lambda
 \qquad{\rm and}\qquad g = \sum_{S,d} b_{S,d} Q_{S,d},$$
then
$$ a_\lambda = \sum_{\{S \mid {\rm type}(S)=\lambda\}} b_{S,d}.$$}

\proof
It is not difficult to see that
it suffices to prove the theorem for the case $g=\xi_\mu$.
For $d$ a positive integer and $S$ a subset of $[d-1]$, define
$$\tilde Q_{S,d} \defeq \sum_{\scriptstyle i_1 \le i_2 \le \cdots \le i_d
  \atop \scriptstyle i_j < i_{j+1}\ {\rm iff}\ j \in S}
  x_{i_1} x_{i_2} \cdots x_{i_d}.$$
Then
$$m_\lambda = \sum_{\{S\mid {\rm type}(S) = \lambda\}} \tilde Q_{S,d}
  \qquad{\rm and}\qquad
  Q_{S,d} = \sum_{T\supseteq S} \tilde Q_{T,d},$$
where in the first summation $d$ is the size of~$\lambda$.
By an inclusion-exclusion argument,
$$m_\lambda = \sum_{\{S\mid {\rm type}(S) = \lambda\}} \;
   \sum_{T\supseteq S} (-1)^{|T| - |S|} Q_{T,d}.$$
Let $q_{\lambda,T,d}$ be the coefficient of~$Q_{T,d}$ in~$m_\lambda$.
We compute
$$ \sum_{\{ T\mid {\rm type}(T) = \nu\}} q_{\lambda,T,d}.$$
Observe that there is a bijection
between subsets of type~$\lambda$
and orderings of the parts of~$\lambda$: given a subset~$S\subseteq [d-1]$
of type~$\lambda$, take the sequence of the lengths of the subwords
of the word $123\ldots d$ obtained by breaking after each
element of~$S$.
Thinking of such subwords as directed paths, we see that for any
fixed $S$ of type~$\lambda$, the number of subsets $T\supseteq S$
such that ${\rm type}(T) = \nu$ is just the number of subsets~$F$
of $E(D_\lambda)$ satisfying $\pi(F)=\nu$,
which from \technical\ is
$${r_\lambda!\over r_\nu!} \,c_{\nu,\lambda}.$$
Now there are
$\ell(\lambda)!/r_\lambda!$ subsets~$S$ of type~$\lambda$,
and if ${\rm type}(S)=\lambda$ and ${\rm type}(T)=\nu$ then
$$(-1)^{|T|-|S|} = ({\rm sgn}\,\nu)({\rm sgn}\,\lambda).$$
Putting all this together, we see that
$$ \sum_{\{ T\mid {\rm type}(T) = \nu\}} \!\!\!q_{\lambda,T,d}
  = {\ell(\lambda)! \over r_\nu!} \,c_{\nu,\lambda} 
    ({\rm sgn}\,\nu)({\rm sgn}\,\lambda).$$
But, again from \technical,
$$\xi_\mu = \sum_\lambda {r_\mu! \over r_\lambda!} \,c_{\lambda,\mu}
   {r_\lambda!\over \ell(\lambda)!} m_\lambda.$$
Hence if $g=\xi_\mu$, then
$$\eqalign{ \sum_{\{S\mid {\rm type}(S) = \nu\}} \!\!\!b_{S,d} &= 
    \sum_\lambda {r_\mu! \over r_\lambda!} \,c_{\lambda,\mu}
    {r_\lambda!\over \ell(\lambda)!} \cdot
    {\ell(\lambda)!\over r_\nu!} \,c_{\nu,\lambda}
    ({\rm sgn}\,\nu)({\rm sgn}\,\lambda)\cr
 &= {r_\mu!\over r_\nu!} \sum_\lambda 
    ({\rm sgn}\,\nu)c_{\nu,\lambda}({\rm sgn}\,\lambda)c_{\lambda,\mu}\cr
 &= \delta_{\mu\nu},\cr}$$
because $\bigl(({\rm sgn}\,\lambda)c_{\lambda,\mu}\bigr)$ is the matrix
of~$\omega$ with respect to the augmented
monomial symmetric function basis (by~[5, Appendix~1, \#9]) and
$\omega$ is an involution.  This completes the proof.\qed

It follows as an immediate corollary that any symmetric function
(such as $\X_G$ or~$s_\lambda$)
that is $Q$-positive is also $\xi$-positive,
and moreover if there is a combinatorial
interpretation of the $Q_{S,d}$-coefficients then it carries
over into a combinatorial interpretation of
the $\xi_\lambda$~coefficients.

We should caution the reader, however,
that $\xi_\lambda$ is \itc{not} $Q$-positive.
Nor is it true that the only $Q_{S,d}$'s in the
$Q_{S,d}$-expansion of~$\xi_\lambda$ with nonzero coefficients
are those with ${\rm type}(S)=\lambda$.
Thus, while \xiquasi\ allows one to translate combinatorial
\itc{interpretations of the coefficients} of the $Q$-expansion
of a symmetric function~$g$ into combinatorial interpretations of
the the coefficients of the $\xi$-expansion of~$g$, there is
no guarantee that combinatorial \itc{proofs} can be so translated.
Some tricky reshuffling of combinatorial information occurs in the
transition from the $Q_{S,d}$'s to the~$\xi_\lambda$'s.
In fact, we do not know of a direct combinatorial proof
that the $\xi_\lambda$-expansion of the Schur functions
enumerates Young tableaux according to descents.

\sec Acknowledgments

This work was supported in part by a National Science Foundation Graduate
Fellowship and a National Science Foundation Postdoctoral Fellowship.
\xiquasi\ first appeared in my M.I.T. doctoral thesis under Richard
Stanley but has not been published before.

\sec References

\item{1.} {F.~Brenti,} Expansions of chromatic polynomials and
log-concavity, \itc{Trans.\ Amer.\ Math.\ Soc.}\ {\bf 332} (1992), 729--756.

\item{2.} {T.~Chow,} The path-cycle symmetric function of a digraph,
\itc{Advances in Math.}\ {\bf 118} (1996), 71--98.

\item{3.} {F.~Chung and R.~Graham,} On the cover polynomial of a
digraph, \itc{J.~Combin.\ Theory (B)} {\bf 65} (1995), 273--290.

\item{4.} {F.~N.~David and M.~G.~Kendall,}
Tables of symmetric functions.~I,
\itc{Biometrika} {\bf 36} (1949), 431--449.

\item{5} {P.~Doubilet,} On the foundations of combinatorial
theory.~VII:
Symmetric functions through the theory of distribution and occupancy,
\itc{Studies in Applied Math.} {\bf 51} (1972), 377--396.

\item{6.} {V.~Gasharov,} Incomparability graphs of $(3+1)$-free
posets are $s$-positive, \itc{Discrete Math.}\ {\bf 157} (1996), 193--197.

\item{7.} {I.~M.~Gessel,}
Multipartite $P$-partitions and inner products of
skew Schur functions, \itc{in} ``Combinatorics and Algebra'' (C.~Greene,
Ed.), Contemporary Mathematics Series, Vol.~34, Amer.\ Math.\ Soc.,
Providence, R.I., 1984, pp.~289--301.

\item{8.} {I.~G.~Macdonald,} \itc{Symmetric Functions and Hall
Polynomials,} Oxford University Press, Oxford, 1979.

\item{9.} {A.~Magid,} Enumeration of convex polyominoes: a generalization
of the Robinson-Schensted correspondence and the dimer problem, Ph.D.
thesis, Brandeis University, 1992.

\item{10.} {B.~E.~Sagan,} \itc{The Symmetric Group: Representations,
Combinatorial Algorithms, and Symmetric Functions,} Wadsworth \&
Brooks/Cole, Pacific Grove, CA, 1991.

\item{11.} {R.~P.~Stanley,} A symmetric function generalization of the
chromatic polynomial of a graph, \itc{Advances in Math.} {\bf 111} (1995),
166--194.

\item{12.} {R.~P.~Stanley,}
\itc{Enumerative Combinatorics,} vol.~1, Wadsworth
\& Brooks/Cole, Pacific Grove, CA, 1986.

\item{13.} {R.~P.~Stanley,} Graph colorings and related symmetric
functions: ideas and applications, preprint.

\item{14.} {T.~S.~Sundquist, D.~G.~Wagner, and J.~West,} A Robinson-Schensted
algorithm for a class of partial orders, \itc{J. Combin.\ Theory Ser.~A}
{\bf 79} (1997), 36--52.
\bye